\documentclass[12pt,a4]{amsart}
\usepackage{amsthm,amsmath,amssymb}
\usepackage{titlesec}
\usepackage{mathrsfs}
\usepackage{comment}
\makeatletter
\def\section{\@startsection{section}{1}%
  \z@{.7\linespacing\@plus\linespacing}{.5\linespacing}%
  {\normalfont\scshape\centering}}
\def\subsection{\@startsection{subsection}{2}%
  \z@{.5\linespacing\@plus.7\linespacing}{-.5em}%
  {\normalfont\bfseries}}
\makeatother
\titleformat*{\section}{\large\bfseries}
\titleformat*{\subsection}{\large\bfseries}
\newtheorem{theorem}{Theorem}[section]

\newtheorem{prop}{Proposition}[section]

\theoremstyle{remark}
\newtheorem{rem}{Remark}
\renewcommand{\Re}{\operatorname{Re}}

\renewcommand{\labelenumi}{{\upshape(\roman{enumi})}}
\makeatletter
\@namedef{subjclassname@2020}{%
  \textup{2020} Mathematics Subject Classification}
\makeatother
\title[the Riemann zeta function over an arithmetic progression]
{A mean-square value of the Riemann zeta function over an arithmetic progression}
\author{Hirotaka Kobayashi}
\date{}
\address{Graduate School of Mathematics, Nagoya University, Furocho, Chikusaku, Nagoya 464-8602, Japan}
\email{m17011z@math.nagoya-u.ac.jp}
\begin{document}

\maketitle

\begin{abstract}
We obtain an asymptotic formula for the second discrete moment of the Riemann zeta function over the arithmetic progression $\frac{1}{2} + in$.
It shows that the first main term is equal to that of the continuous mean value.
\end{abstract}

\section{Introduction}
In this paper, we shall consider the mean-square
\begin{equation*}
\sum_{n \leq T} \left|\zeta \left(\frac{1}{2} + in \right) \right|^2.
\end{equation*}
This study is one of the attempts to enrich our knowledge about the vertical distribution of the Riemann zeta function $\zeta(s)$.
It was Putnam \cite{Put} who first cultivated the study of the values of $\zeta(s)$ on vertical arithmetic progressions.
He showed that there is no infinite arithmetic progression of non-trivial zeros of $\zeta(s)$.
In this direction, there is an important conjecture called the linear independence conjecture, which states that the ordinates of non-trivial zeros of $\zeta(s)$ are linearly independent over $\mathbb{Q}$. In 1942, Ingham \cite{Ing} found out the relation between the linear independence conjecture and the oscillations of $M(x) = \sum_{n \leq x} \mu(n)$, where $\mu(n)$ is the M\"{o}bius function. He showed that the linear independent conjecture implies the failure of the inequality $M(x) \ll x^{1/2}$. For this reason, many mathematicians doubt the inequality.

We may need much more progression to solve the problem. However, we have an easier conjecture in this direction, that is, there are no non-trivial zeros of $\zeta(s)$ in arithmetic progression. To consider this problem, discrete moments play an important role. Martin and Ng \cite{Ma&Ng} attacked this conjecture for Dirichlet $L$-functions by considering some kinds of discrete means of Dirichlet $L$-functions. Later, by a similar approach,  Li and Radziwi\l \l \cite{Li&Ra} showed that at least one third of the values on arithmetic progression does not vanish. One of their results (Theorem 1) stated that for fixed arbitrary $A > 0$, as $T \to \infty$,
\begin{equation*}
\begin{split}
&\quad
\sum_{n} \left|\zeta \left(\frac{1}{2} + in \right) B\left(\frac{1}{2} + in \right) \right|^2 \cdot \phi \left(\frac{n}{T} \right) \\
&=
\int_{\mathbb{R}} \left|\zeta \left(\frac{1}{2} + it \right) B\left(\frac{1}{2} + it \right) \right|^2 \cdot \phi \left(\frac{t}{T} \right) dt + O_{A}(T (\log T)^{-A}),
\end{split}
\end{equation*}
where $B(s)$ is a Dirichlet polynomial and $\phi(\cdot)$ is a smooth compactly supported function.
This clarify the notable correspondence of the discrete mean to the continuous one.

\"{O}zbek and Steuding \cite{Oz&St2} proved asymptotic formulas for the first discrete moment of $\zeta(s)$ on certain vertical arithmetic progressions inside the critical strip. The first discrete moment have been studied recently by \"{O}zbek, Steuding and Wegert (see \cite{St&We} and \cite{Oz&St1}). Especially, in \cite{Oz&St1}, they proved showed that
\begin{equation*}
\lim_{T \to \infty} \frac{1}{T}\sum_{0 \leq n < T} \zeta(s_0 + in\delta)
=
\begin{cases}
(1 - l^{-s_0})^{-1} & \text{if} \ \delta = \frac{2\pi q}{\log l}, q \in \mathbb{N}, 2 \leq l \in \mathbb{N}, \\
1 & \text{otherwise},
\end{cases}
\end{equation*}
where $s_0$ may be any complex number with real part in $(0,1)$.

The object of this paper is to prove an asymptotic formula for the discrete mean-square of the Riemann zeta function on a vertical arithmetic progression.
\begin{theorem}
We have, as $T \to \infty$,
\begin{equation}\label{th}
\sum_{n \leq T} \left|\zeta \left(\frac{1}{2} + in \right) \right|^2 = T\log \frac{T}{2\pi} + O \left(T \log T \exp \left(-C \frac{\log_2 T}{\log_3 T}\right) \right)
\end{equation}
for some constant $C > 0$.
\end{theorem}

This should be compared with the continuous mean-square
\begin{equation*}
\int_{0}^{T} \left|\zeta \left(\frac{1}{2} + it \right) \right|^2 dt = T \log \frac{T}{2\pi} + (2\gamma - 1)T + E(T),
\end{equation*}
where $\gamma$ is Euler's constant and $E(T)$ is an error term. This reveals that the discrete mean value (\ref{th}) equals to the continuous mean value as $T \to \infty$ asymptotically. Moreover, if a conjecture on the irrationality measure of $e^{2\pi m}$ is true, we can see that the second main term also coincides.

The error term estimate in (\ref{th}) is due to Bugeaud and Ivi\'{c} \cite{Bug&Iv}. They studied the discrete mean value of $E(T)$. In their proof, the sum
\begin{equation*}
\sum_{m \leq Te^{-2\pi k}}d(m) \exp(-2\pi mi e^{2\pi k})
\end{equation*}
appears, and they estimate the upper bound. Similarly, in our proof, we have to estimate the same sum. Thus, the same problem arises in the discrete mean value of $E(T)$ and the discrete mean-square of $\zeta(s)$. The estimate requires bounds for the irrationality measure of $e^{2\pi k}$ and for the partial quotients in its continued fraction expansion.

Our starting point is the approximate functional equation of $\zeta^2(s)$:
\begin{equation}\label{app-zeta2}
\zeta^2(s) = \sum_{n \leq \frac{t}{2\pi}} \frac{d(n)}{n^s} + \chi^2(s) \sum_{n \leq \frac{t}{2\pi}} \frac{d(n)}{n^{1-s}} + R \left(s ; \frac{t}{2\pi} \right),
\end{equation}
where $\chi(s) = 2^s \pi^{s-1} \sin(\pi s/2)\Gamma(1-s)$ and $R(s ; t/2\pi)$ is the error term.
Motohashi \cite{Mo1}, \cite{Mo2} proved that
\begin{equation}\label{R-d-simple}
\chi(1-s) R \left(s ; \frac{t}{2\pi} \right) = -\sqrt{2}\left( \frac{t}{2\pi} \right)^{-\frac{1}{2}}\Delta \left( \frac{t}{2\pi} \right) + O(t^{-1/4}),
\end{equation}
where $\Delta(t/2\pi)$ is the error term in the Dirichlet divisor problem, defined by
\begin{equation*}
\Delta(x) = \sideset{}{'} \sum_{n \leq x} d(n) - x(\log x + 2\gamma - 1) - \frac{1}{4}.
\end{equation*}
Here $\sum'$ indicates that the last term is to be halved if $x$ is an integer. We note that Jutila \cite{Ju} gave another proof of Motohashi's result (\ref{R-d-simple}).

\begin{rem}
van Frankenhuijsen \cite{van1} gave an explicit bound for the length of arithmetic progressions of non-trivial zeros of the Riemann zeta function. Later, he \cite{van2} improved the bound.
\end{rem}

\begin{rem}
Good \cite{Go} proved asymptotic formulas for the fourth moments of the Riemann zeta function on arbitrary arithmetic progressions to the right of the critical line. For example, he showed that for $\sigma>1/2$
\begin{equation*}
\sum_{0 \leq n < T} |\zeta(\sigma + i nd)|^4 = T \sum_{m = 1}^{\infty} \frac{d(m)^2}{m^{2\sigma}} + o(T) \quad (T \to \infty),
\end{equation*}
 where $d$ is not of the form $2\pi l (\log k_1/k_2)^{-1}$ with integral $l \neq 0$ and positive integers $k_1 \neq k_2$.
\end{rem}

\section{the proof of theorem}

By (\ref{app-zeta2}), (\ref{R-d-simple}) and the functional equation $\zeta(1-s) = \chi(1-s) \zeta(s)$, we have
\begin{equation*}
\begin{split}
\zeta(s)\zeta(1-s)
&=
\chi(1-s) \sum_{n \leq \frac{t}{2\pi}} \frac{d(n)}{n^s} + \chi (s) \sum_{n \leq \frac{t}{2\pi}} \frac{d(n)}{n^{1-s}} \\
&\quad
-\sqrt{2}\left( \frac{t}{2\pi} \right)^{-\frac{1}{2}}\Delta \left( \frac{t}{2\pi} \right) + O(t^{-1/4}).
 \end{split}
\end{equation*}

It is known that $\Delta(t) \ll t^{1/3}$. Thus, taking $s = 1/2 + it$, we have
\begin{equation*}
\left|\zeta \left(\frac{1}{2} + it \right) \right|^2 = 2\Re \chi \left(\frac{1}{2} - it \right) \sum_{n \leq \frac{t}{2\pi}} \frac{d(n)}{n^{1/2 + it}} + O(t^{-1/6}).
\end{equation*}
Hence we consider
\begin{equation*}
\begin{split}
\sum_{T < n \leq 2T}\left|\zeta \left(\frac{1}{2} + in \right) \right|^2
&=
2 \Re \sum_{T < n \leq 2T} \chi \left(\frac{1}{2} - in \right) \sum_{m \leq \frac{n}{2\pi}} \frac{d(m)}{m^{1/2 + in}} \\
&\quad
+ O \left(\sum_{T < n \leq 2T}n^{-1/6}\right).
\end{split}
\end{equation*}

As for the error term, it is clear that
\begin{equation*}
\sum_{T < n \leq 2T} n^{-1/6} \ll T^{5/6}.
\end{equation*}

To consider the main term, we note the following fact:
\begin{equation*}
\chi(1-s)=e^{-\pi i/4} \left(\frac{t}{2\pi} \right)^{\sigma - 1/2}\exp \left(it \log \frac{t}{2\pi e} \right)(1+ O (t^{-1}))
\end{equation*}
for $\sigma$ fixed and $t \geq 1$. Using this, we have
\begin{equation*}
\begin{split}
&\quad
\sum_{T < n \leq 2T} \chi \left(\frac{1}{2} - in \right) \sum_{m \leq \frac{n}{2\pi}} \frac{d(m)}{m^{1/2 + in}} \\
&=
e^{-\pi i/4}\sum_{T < n \leq 2T} \exp \left(in \log \frac{n}{2\pi e} \right)(1+ O (n^{-1})) \sum_{m \leq \frac{n}{2\pi}} \frac{d(m)}{m^{1/2 + in}} \\
&=
e^{-\pi i/4} \sum_{m \leq \frac{T}{\pi}} \frac{d(m)}{m^{1/2}} \sum_{\max(2\pi m, T) < n \leq 2T} \exp \left(in \log \frac{n}{2\pi em} \right) \\
&\quad
+ O(T^{1/2} \log T) \\
&=
e^{-\pi i/4} \sum_{ m \leq \frac{T}{2\pi}} \frac{d(m)}{m^{1/2}} \sum_{T < n \leq 2T} \exp \left(in \log \frac{n}{2\pi em} \right) \\
&\quad
+ e^{-\pi i/4} \sum_{\frac{T}{2\pi} < m \leq \frac{T}{\pi}} \frac{d(m)}{m^{1/2}} \sum_{2\pi m < n \leq 2T} \exp \left(in \log \frac{n}{2\pi em} \right) \\
&\quad
+ O(T^{1/2} \log T) \\
&=
 e^{-\pi i/4}(S_1 + S_2) + O(T^{1/2} \log T),
\end{split}
\end{equation*}
say. To obtain the second equality, we note that
\begin{equation*}
\sum_{m \leq T}\frac{d(m)}{m^{1/2}} \ll T^{1/2}\log T.
\end{equation*}

For the convenience, we put
\begin{equation*}
f(x) := \frac{x}{2\pi} \log \frac{x}{2\pi em}
\end{equation*}
and
\begin{equation*}
g_k(x) := f(x) - kx.
\end{equation*}
To calculate $S_1$ and $S_2$, we use the saddle-point method.

\subsection{An estimate of $S_1$}
Since derivatives of $f(x)$ are
\begin{equation*}
f'(x) = \frac{1}{2\pi} \log \frac{x}{2\pi m}, \quad f''(x) = \frac{1}{2\pi x},
\end{equation*}
we have (see Proposition 8.7. in \cite{Iw&Ko})
\begin{equation*}
\sum_{T < n \leq 2T}e^{2\pi if(n)} = \sum_{\frac{1}{2\pi}\log \frac{T}{2\pi m} - \theta < k < \frac{1}{2\pi}\log \frac{T}{\pi m} + \theta} \int_{T}^{2T} e(g_k(x)) dx + O(1),
\end{equation*}
where $\theta$ is any sufficiently small positive number and $e(x) := \exp(2\pi ix)$. Then $S_1$ can be rewritten as
\begin{equation*}
\begin{split}
&\quad
\sum_{m \leq \frac{T}{2\pi}} \frac{d(m)}{m^{1/2}} \sum_{\frac{1}{2\pi}\log \frac{T}{2\pi m} - \theta < k < \frac{1}{2\pi}\log \frac{T}{\pi m} + \theta} \int_{T}^{2T} e(g_k(x)) dx +O(T^{1/2}\log T) \\
&=
\sum_{0 \leq k \leq \frac{1}{2\pi}\log \frac{T}{\pi} + \theta} \sum_{\frac{T}{2\pi}e^{-2\pi(k + \theta)} < m < \frac{T}{\pi}e^{-2\pi(k - \theta)}} \frac{d(m)}{m^{1/2}} \int_{T}^{2T} e(g_k(x)) dx \\
&\quad
+ O(T^{1/2}\log T).
\end{split}
\end{equation*}
We note that the saddle-point of $g_k(x)$, $2\pi me^{2\pi k}$, is in $(Te^{-2\pi \theta}, 2Te^{2\pi \theta})$ by the condition in the inner sum.
We divide $(\frac{T}{2\pi}e^{-2\pi(k + \theta)}, \frac{T}{\pi}e^{-2\pi(k - \theta)})$ into the following five intervals:
\begin{enumerate}
\setlength{\itemsep}{0.5cm}
\item
$\left(\frac{T}{2\pi}e^{-2\pi(k + \theta)}, \frac{T}{2\pi}e^{-2\pi k} - 1 \right]$,
\item
$\left(\frac{T}{2\pi}e^{-2\pi k} - 1, \frac{T}{2\pi}e^{-2\pi k} + 1 \right)$,
\item
$\left[\frac{T}{2\pi}e^{-2\pi k} + 1, \frac{T}{\pi}e^{-2\pi k} - 1 \right]$,
\item
$\left(\frac{T}{\pi}e^{-2\pi k} - 1, \frac{T}{\pi}e^{-2\pi k} + 1 \right)$,
\item
$\left[\frac{T}{\pi}e^{-2\pi k} + 1, \frac{T}{\pi}e^{- 2\pi (k - \theta)} \right)$.
\end{enumerate}

(i)

By the first derivative test, we have
\begin{equation*}
\left|\int_{T}^{2T}e(g_k(x))dx \right| \leq \frac{8\pi}{\log \frac{T}{2\pi me^{2\pi k}}},
\end{equation*}
and
\begin{equation*}
\left|\log \frac{T}{2\pi me^{2\pi k}} \right| = \left|- \log \left(1 - \frac{T - 2\pi me^{2\pi k}}{T}\right) \right| \sim \frac{T - 2\pi me^{2\pi k}}{T}.
\end{equation*}
Therefore, in this case, the contribution is
\begin{equation*}
\begin{split}
&\ll
\sum_{0 \leq k \leq \frac{1}{2\pi}\log \frac{T}{\pi}} \sum_{m \leq \frac{T}{2\pi}e^{-2\pi k} - 1} \frac{d(m)}{m^{1/2}} \frac{T}{T - 2\pi me^{2\pi k}} \\
&\ll
T^{\varepsilon} \sum_{0 \leq k \leq \frac{1}{2\pi}\log \frac{T}{\pi}} Te^{-2\pi k} \sum_{m \leq \frac{T}{2\pi}e^{-2\pi k} - 1} \frac{m^{1/2}}{m(\frac{T}{2\pi}e^{-2\pi k} - m)} \\
&\ll
T^{\varepsilon} \sum_{0 \leq k \leq \frac{1}{2\pi}\log \frac{T}{\pi}} \sum_{m \leq \frac{T}{2\pi}e^{-2\pi k} - 1} m^{1/2}\left(\frac{1}{m} + \frac{1}{\frac{T}{2\pi}e^{-2\pi k} - m} \right) \\
&\ll
T^{1/2 + \varepsilon} \sum_{0 \leq k \leq \frac{1}{2\pi}\log \frac{T}{\pi}} e^{-\pi k} \sum_{m \leq \frac{T}{2\pi}}\frac{1}{m} \ll T^{1/2 + \varepsilon}.
\end{split}
\end{equation*}

(v)

In a similar manner, we can see that the contribution of this case is $\ll T^{1/2 + \varepsilon}$.

(ii), (iv)

First we note that the number of $m$ in each intervals are at most two and
\begin{equation*}
m \asymp Te^{-2\pi k} \ll T.
\end{equation*}
By the second derivative test, we have
\begin{equation*}
\left|\int_{T}^{2T}e(g_k(x))dx \right| \leq 16\sqrt{\pi T}.
\end{equation*}
Hence we can see that
\begin{equation*}
T^{1/2} \sum_{0 \leq k \leq \frac{1}{2\pi}\log \frac{T}{\pi} + \theta} \sum_{\substack{\frac{T}{2\pi}e^{-2\pi k} - 1 < m < \frac{T}{2\pi}e^{-2\pi k} + 1 \\ \frac{T}{\pi}e^{-2\pi k} - 1 < m < \frac{T}{\pi}e^{-2\pi k} + 1}} \frac{d(m)}{m^{1/2}} \ll T^{1/2+\varepsilon}.
\end{equation*}

(iii)

In this case, the sum of $k$ start from $1$. If $k = 0$, then we have that $\frac{T}{2\pi} + 1 < m$. However, this is impossible, because we consider the case $m \leq \frac{T}{2\pi}$.
Using the saddle-point method (see for example Corollary 8.15. in \cite{Iw&Ko}), we have
\begin{equation*}
\begin{split}
\int_{T}^{2T}e(g_k(x))dx 
&=
e^{\pi i/4} e^{-2\pi mi e^{2\pi k}} 2\pi e^{\pi k} \sqrt{m} \\
&\quad
+O \left(\frac{T}{2T - 2\pi me^{2\pi k}} + \frac{T}{2\pi me^{2\pi k} - T} + 1 \right).
\end{split}
\end{equation*}

By the same argument with the case (i) and (v), we see that the contribution of the error term is $\ll T^{1/2 + \varepsilon}$. We have to consider
\begin{equation}\label{irrational-sum}
\sum_{1 \leq k < \frac{1}{2\pi}\log \frac{T}{\pi}} e^{\pi k} \sum_{\frac{T}{2\pi}e^{-2\pi k} + 1 < m < \frac{T}{\pi}e^{-2\pi k} - 1} d(m) \exp(-2\pi mi e^{2\pi k}).
\end{equation}
However, Bugeaud and Ivi\'{c} \cite{Bug&Iv} have already studied a quite similar sum (\cite[(2.6)]{Bug&Iv}).
They applied the following two results:
\begin{enumerate}
\renewcommand{\labelenumi}{{\upshape(\alph{enumi})}}
\item
The functional equation by Wilton \cite{Wil}
\begin{equation*}
D(x, \eta) = \eta^{-1} D(\eta^2 x, -\eta^{-1}) + O(x^{1/2} \log x), \quad D(x, \eta) = \sum_{m \leq x} d(m) e^{2\pi i\eta m},
\end{equation*}
where $\eta^2 x\ll 1$ and $0 < \eta \leq 1$ is real.

\item
A bound for the irrationality measure of $e^{\pi k}$ by Waldschmidt \cite[p.473]{Wal}

For any positive integers $m$, $p$ and $q$ with $p \leq 3$, we have
\begin{equation*}
\left|e^{\pi k} - \frac{p}{q} \right| > \exp \{ -2^{72}(\log 2k)(\log p)(\log \log p) \}.
\end{equation*}
\end{enumerate}

Now let
\begin{equation*}
e^{\pi k} = [a_0(k); a_1(k), \dots]
\end{equation*}
be the expansion of $e^{\pi k}$ as a continued fraction for any non-zero integer $k$.
From (b), they proved
\begin{prop}[Lemma 1 in \cite{Bug&Iv}]
There exists an absolute positive constant $c$ such that
\begin{equation*}
\log \log a_n(k) < c(n + \log |k|) \log (n + \log |k|)
\end{equation*}
holds for any non-zero integer $k$ and any non-negative integer $n$.
\end{prop}
Using these results, they showed that
\begin{equation*}
\sum_{k \ll \log T} \frac{e^{\pi k}}{k} \sum_{m \leq T e^{-2\pi k}} d(m) \exp(-2\pi mi e^{2\pi k}) \ll T\log T \exp \left(-C \frac{\log_2 T}{\log_3 T}\right),
\end{equation*}
for some constant $C > 0$. In our sum (\ref{irrational-sum}), the factor $1/k$ does not appear, but this does not influence the estimate. Thus we can bound the sum (\ref{irrational-sum}) by
\begin{equation*}
T\log T \exp \left(-C \frac{\log_2 T}{\log_3 T}\right)
\end{equation*}
for some constant $C > 0$. Consequently, $S_1$ is bounded by
\begin{equation}\label{error}
T\log T \exp \left(-C \frac{\log_2 T}{\log_3 T}\right).
\end{equation}

\begin{rem}
As Bugeaud and Ivi\'{c} \cite{Bug&Iv} referred, from the result of Wilton \cite{Wil}, if $a_n(k)$ satisfy $a_n(k) \ll n^{1 + K} \ (K \geq 0)$, then
\begin{equation*}
\sum_{m \leq x}d(m) \exp(2\pi i m e^{2\pi k}) \ll x^{1/2} \log^{2 + K} x.
\end{equation*}
If this estimate is verified, we can obtain the second main term $(2\gamma - 1)T$ in (\ref{th}).
\end{rem}

\subsection{A contribution of $S_2$}
As in the case $S_1$, we have
\begin{equation*}
\sum_{2\pi m < n \leq 2T}e^{2\pi if(n)} = \sum_{0 \leq k < \frac{1}{2\pi}\log \frac{T}{\pi m} + \theta} \int_{2\pi m}^{2T} e(g_k(x)) dx + O(\theta^{-1} + \log \log T),
\end{equation*}
where $\theta$ is any sufficiently small positive number.

Since
\begin{equation*}
\sum_{m \leq T} \frac{d(m)}{m^{1/2}} \ll T^{1/2}\log T,
\end{equation*}
the contribution of the error term is $\ll T^{1/2} \log T \log \log T$.

Finally, we calculate
\begin{equation*}
\sum_{\frac{T}{2\pi} < m \leq \frac{T}{\pi}} \frac{d(m)}{m^{1/2}} \sum_{0 \leq k < \frac{1}{2\pi}\log \frac{T}{\pi m} + \theta} \int_{2\pi m}^{2T} e(g_k(x)) dx.
\end{equation*}
However, the conditions for $m$ and $k$ imply that $k = 0$. Therefore the above is
\begin{equation*}
\sum_{\frac{T}{2\pi} < m \leq \frac{T}{\pi}} \frac{d(m)}{m^{1/2}} \int_{2\pi m}^{2T} e(f(x)) dx.
\end{equation*}

The saddle-point of $f(x)$ is $2\pi m$, thus, the saddle-point method leads to
\begin{equation*}
\int_{2\pi m}^{2T} e(f(x)) dx = e^{\pi i/4} \pi \sqrt{m} + O \left(\frac{T}{2T - 2\pi m} + 1 \right).
\end{equation*}
Therefore 
\begin{equation}\label{main}
\sum_{\frac{T}{2\pi} < m \leq \frac{T}{\pi}} \frac{d(m)}{m^{1/2}} \int_{2\pi m}^{2T} e(f(x)) dx = e^{\pi i/4} \pi \sum_{\frac{T}{2\pi} < m \leq \frac{T}{\pi}} d(m) + O(T^{1/2 + \varepsilon}).
\end{equation}

\subsection{conclusion}
By (\ref{error}) and (\ref{main}), we can see that
\begin{equation*}
\begin{split}
&\quad
\sum_{T < n \leq 2T} \chi \left(\frac{1}{2} - in \right) \sum_{m \leq \frac{n}{2\pi}} \frac{d(m)}{m^{1/2 + in}} \\
&=
\pi \sum_{\frac{T}{2\pi} < m \leq \frac{T}{\pi}} d(m) + O \left(T \log T \exp \left(-C \frac{\log_2 T}{\log_3 T}\right) \right),
\end{split}
\end{equation*}
and so,
\begin{equation*}
\begin{split}
&\quad
\sum_{T < n \leq 2T}\left|\zeta \left(\frac{1}{2} + in \right) \right|^2 \\
&=
2 \pi \sum_{\frac{T}{2\pi} < m \leq \frac{T}{\pi}} d(m) + O \left(T \log T \exp \left(-C \frac{\log_2 T}{\log_3 T}\right) \right).
\end{split}
\end{equation*}

Finally, replacing $T$ by $\frac{T}{2}, \frac{T}{4}$, and so on, and adding we have
\begin{equation*}
\begin{split}
\sum_{n \leq T}\left|\zeta \left(\frac{1}{2} + in \right) \right|^2
&=
2 \pi \sum_{ m \leq \frac{T}{2\pi}} d(m) + O \left(T \log T \exp \left(-C \frac{\log_2 T}{\log_3 T}\right) \right) \\
&=
T\log \frac{T}{2\pi} + + O \left(T \log T \exp \left(-C \frac{\log_2 T}{\log_3 T}\right) \right).
\end{split}
\end{equation*}

This completes the proof.

\end{document}